\documentclass[11pt]{article}
\usepackage{amsmath, amssymb, amsfonts, amsthm, amscd, vmargin}
\usepackage[latin1]{inputenc}
\usepackage[all]{xy}
\usepackage{graphicx}
\usepackage{enumerate}

\setmarginsrb{2.4cm}{4cm}{2.4cm}{3.3cm}{0cm}{0cm}{0cm}{1.5cm}

\theoremstyle{plain}
\newtheorem{teo}{Theorem}[section]
\newtheorem{lem}[teo]{Lemma}
\newtheorem{prop}[teo]{Proposition}

\theoremstyle{definition}
\newtheorem{defi}[teo]{Definition}

\newtheorem{rem}[teo]{Remark}

\newenvironment{preu}{{\flushleft \it Proof of Theorem
    \ref{eclatement}.}}{\hfill
$\square$ \vspace{2mm}}
\numberwithin{equation}{teo}
\newcommand{\Aut}{\operatorname{Aut}}

\renewcommand{\Im}{\operatorname{Im}}

\newcommand{\Cbb}{{\mathbb C}}
\newcommand{\Qbb}{{\mathbb Q}}
\newcommand{\Zbb}{{\mathbb Z}}

\newcommand{\Pbb}{{\mathbb P}}

\newcommand{\Obb}{{\mathbb O}}

\newcommand{\lra}{\longrightarrow}

\newcommand{\Zt}{{\widetilde{Z}}}
\newcommand{\Xt}{{\widetilde{X}}}

\def\cO{{\cal O}}

\begin{document}

\title{Local rigidity
of quasi-regular varieties}

\author{Boris Pasquier and Nicolas Perrin}
\maketitle

\begin{abstract}
For a $G$-variety $X$ with an open orbit, we define its boundary
$\partial X$ as the complement of the open orbit. The action sheaf
$S_X$ is the subsheaf of the tangent sheaf made of vector fields
tangent to $\partial X$. We prove, for a large family of smooth
spherical varieties, the vanishing of the cohomology groups
$H^i(X,S_X)$ for $i>0$, extending results of F. Bien and M. Brion
\cite{BB96}.

We apply these results to study the local rigidity of the smooth
projective varieties with Picard number one classified in
\cite{Pa08}.
\end{abstract}



\tableofcontents

\section*{Introduction}

Let $X$ be a complex algebraic variety. Denote by $T_X$ its tangent
bundle. If $X$ is a flag variety, it is well-known that
$H^i(X,T_X)=0$ for any $i\geq 1$ (see \cite{demazure}). By
Kodaira-Spencer theory, the vanishing of $H^1(X,T_X)$ implies that
$X$ is locally rigid, {\it i.e.} admits no local deformation of its
complex structure.

Let $G$ be a connected reductive algebraic group over $\Cbb$ and $X$
be a smooth projective $G$-variety. Then $X$ is said to be regular
if it is smooth, spherical ({\it i.e.} has an open orbit under the
action of a Borel subgroup $B$ of $G$) without color ({\it i.e.}
every irreducible $B$-stable divisor containing a $G$-orbit is
$G$-stable). F.~Bien and M.~Brion proved in \cite{BB96} that regular
varieties are not locally rigid in general but they have a weaker
rigidity property. Indeed, let $X$ be a spherical variety, we denote
by $\Omega$ its open $G$-orbit and by $\partial X$ its complement
(we shall call it the boundary of $X$). Denote by $S_X$ the action
sheaf of $X$ {\it i.e.} the subsheaf of $T_X$ made of vector fields
tangent to $\partial X$. Then, combining Theorem 4.1 in \cite{Kn94}
and Proposition 2.5 in \cite{BB96}, the following result holds:

\begin{teo}\label{BBK}
Let $X$ be a projective regular variety, then $H^i(X,S_X)=0$ for any
$i>0$.
\end{teo}

In this paper, we generalise this result for a larger family of
spherical varieties.

\begin{defi}
Let $X$ be a smooth projective spherical variety. Denote by $n$ its
rank, {\it i.e.} the minimal codimension of $U$-orbits in the open
$G$-orbit of $X$ (where $U$ is the unipotent radical of the Borel
subgroup $B$).
The variety $X$ is said to be quasi-regular if the following
conditions hold:
\begin{itemize}
\item[(QR1)] any irreducible component of the boundary $\partial X$
is a smooth (spherical) variety of rank $n-1$;
\item[(QR2)] any $G$-orbit closure $Y$ is the intersection of
the $Z_i$ containing $Y$.
\end{itemize}
\end{defi}

\begin{rem}
A regular variety is quasi-regular (see \cite{BB96}). Remark also
that our conditions are direct generalisations of the conditions
defining regular varieties in \cite{BB96} (the fact that the
boundary components intersect transversally follows from condition
(QR1), see Lemma \ref{transverse}).
\end{rem}

We will see (Lemma \ref{horospherical}) that the family of
quasi-regular varieties contains all smooth horospherical varieties
({\it i.e.} spherical varieties such that $\Omega$ is a torus bundle
over a flag variety, for more details see also \cite{PaSMF}) and all
smooth spherical varieties of rank~1.
%
The main result of the paper is the following:

\begin{teo}\label{eclatement}
Let $X$ be a quasi-regular spherical variety, then $H^i(X,S_X)=0$ for any $i>0$.
\end{teo}

%

Apart from generalising Theorem \ref{BBK}, a motivation for Theorem
\ref{eclatement} is to prove that certain spherical varieties of
rank~1 are indeed locally rigid. In \cite{Pa08}, the classification
of horospherical varieties of Picard number one was achieved. More
generally, in \emph{loc. cit.}, all smooth projective two-orbits
varieties with Picard number one that still have two orbits after
blowing-up the closed orbit were classified. We shall say that such
a variety satisfies the condition $(\dag)$. The varieties satisfying
$(\dag)$ are spherical of rank~1. Let us describe them more
precisely.

The horospherical varieties satisfying $(\dag)$ have three
$G$-orbits, the open orbit $\Omega$ and two closed orbits $Y\simeq
G/P_Y$ and $Z\simeq G/P_Z$ (but they have only two orbits under the
action of their automorphism group). These varieties are classified
by the triples $(G,P_Y,P_Z)$ in the following list (see
\cite[Th.0.1]{Pa08}, we take the notation of \cite{Bo75} for
fundamental weights $\varpi_i$ and $P(\varpi_i)$ is the associated
parabolic subgroup):
\begin{enumerate}
\item $(B_m,P(\varpi_{m-1}),P(\varpi_m))$ with $m\geq 3$
\item $(B_3,P(\varpi_1),P(\varpi_3))$
\item $(C_m,P(\varpi_i),P(\varpi_{i+1}))$ with $m\geq 2$ and
$i\in\{1,\ldots,m-1\}$
\item $(F_4,P(\varpi_2),P(\varpi_3))$
\item $(G_2,P(\varpi_2),P(\varpi_1))$
\end{enumerate}
We denote by $X^1(m)$, $X^2$, $X^3(m,i)$, $X^4$ and $X^5$ the corresponding
varieties.

There are only two non horospherical varieties satisfying $(\dag)$.
We denote them by $\mathbb{X}_1$ resp. $\mathbb{X}_2$. The
corresponding group $G$ is $F_4$ resp.
$G_2\times\operatorname{PSL(2)}$ (and they have only two
$G$-orbits). See \cite[Definitions 2.11 and 2.12]{Pa08} for explicit
definitions of these varieties.

The varieties $X_3(m,i)$ are the odd symplectic grassmannians
studied by I.~Mihai \cite{Mi05}. They have many nice geometric
properties and are in particular locally rigid. It is thus natural
to ask if the other varieties satisfying $(\dag)$ are also locally
rigid. The answer is as follows:

\begin{teo}\label{deformations}
Assume that $X$ satisfy $(\dag)$, then we have the alternative:
\begin{itemize}\item if $X=X^5$, then $H^1(X,T_X)=\Cbb$ and $H^i(X,T_X)=0$
for any $i\geq 2$;
\item if $X\neq X^5$, then $H^i(X,T_X)=0$ for any $i\geq 1$.
\end{itemize}
\end{teo}

We shall also prove that the non trivial local deformation of the
horospherical $G_2$-variety $X^5$ comes from an actual deformation
to a variety homogeneous under $G_2$ (see Proposition
\ref{deformation}).

\vskip 0.2 cm

Finally we prove (see Proposition \ref{aut}) a characterisation of
homogeneity for spherical varieties of rank~1. This simplifies some
of the proofs given in \cite{Pa08}.
%

\section{Cohomology of the action sheaf}

In this section we will prove Theorem \ref{eclatement}, using
Theorem \ref{BBK}. We will relate any quasi-regular variety to a
regular variety by blow-ups of irreducible components of the
boundary $\partial X$ of $X$.
We first need to recall some notation
and basic facts on spherical varieties.

\subsection{Spherical varieties}

Let $G$ be a reductive connected algebraic group. Let $H$ be a
closed subgroup of $G$. The homogeneous space $G/H$ is said to be
spherical if it has an open orbit under the action of a Borel
subgroup $B$ of $G$. A $G/H$-embedding is a normal $G$-variety that
contains an open orbit isomorphic to $G/H$. Then the spherical
varieties are the $G/H$-embeddings with $G/H$ spherical.

Now, fix a spherical homogeneous space $G/H$ of rank $n$. Then $G/H$-embeddings
have been classified in terms of colored fans by D.~Luna and T.~Vust
\cite{LV83}. Let us recall a part of this theory, see \cite{Kn91}
and \cite{Br97} for more details.

We denote by $\mathcal{D}$ the set of irreducible $B$-stable divisors
of $G/H$. An element of $\mathcal{D}$ is called a color. Let $M$ be
the lattice of all characters $\chi$ of $B$ such that there exists a
non-zero element $f\in\Cbb(G/H)$ such that for all $b\in B$ and $x\in
G/H$ we have $f(bx)=\chi(b)f(x)$ (remark that such a $f$ is unique up
to a scalar). Denote by $N$ the dual lattice of $M$ and let
$N_\Qbb=N\otimes_\Zbb\Qbb$. Note that $N$ and $M$ are lattices of rank $n$.

Let $D\in \mathcal{D}$, then the associated $B$-stable valuation
$v_D:\Cbb(G/H)\to\Zbb$ defines an element of $N_\Qbb$, denoted by
$\rho(D)$.
Denote by $\mathcal{V}$ the image in $N$ of the set of $G$-stable
valuations of $X$. It is a polyhedral convex cone of $N_\Qbb$.

\begin{defi}
(\i) A colored cone is a pair $(\mathcal{C},\mathcal{F})$ with
$\mathcal{C}\subset N_\Qbb$ and $\mathcal{F}\subset\mathcal{D}$ having
the following properties:
\begin{itemize}
\item $\mathcal{C}$ is a convex cone generated by $\rho(\mathcal{F})$
  and finitely many elements of $\mathcal{V}$;
\item the relative interior of $\mathcal{C}$ intersects $\mathcal{V}$
  non trivially;

\item $\mathcal{C}$ contains no lines and $0\not\in\rho(\mathcal{F})$.
\end{itemize}

(\i\i) A colored face of a colored cone $(\mathcal{C},\mathcal{F})$
is a pair $(\mathcal{C'},\mathcal{F'})$ such that $\mathcal{C'}$ is
a face of $\mathcal{C}$, the relative interior of $\mathcal{C'}$
intersects non trivially $\mathcal{V}$ and $\mathcal{F'}$ is the
subset of $\mathcal{F}$ of elements $D$ satisfying
$\rho(D)\in\mathcal{C'}$.

(\i\i\i) A colored fan is a finite set $\mathbb{F}$ of colored cones
with the following properties:
\begin{itemize}
\item every colored face of a colored cone of $\mathbb{F}$ is in $\mathbb{F}$;
\item for all $v\in\mathcal{V}$, there exists at most one
  $(\mathcal{C},\mathcal{F})\in\mathbb{F}$ such that $v$ is in the
  relative interior of $\mathcal{C}$.
\end{itemize}

(\i v) The support of a colored fan $\mathbb{F}$ is the set of
elements of $\mathcal{V}$ contained in the cone of a colored cone of
$\mathbb{F}$.
A color of a colored cone of $\mathbb{F}$ is an element of
$\mathcal{D}$ such that there exists
$(\mathcal{C},\mathcal{F})\in\mathbb{F}$ such that $D\in\mathcal{F}$.
\end{defi}

\begin{teo}[Luna-Vust]\label{Luna-Vust}
There is a bijection $X\longmapsto\mathbb{F}(X)$ between the set of
isomorphism classes of $G/H$-embeddings and the set of colored fans.
Furthermore, under this isomorphism, we have the following
properties:
\begin{itemize}
\item Let $X$ be a $G/H$-embedding.

There exists a bijection between the set of $G$-orbits of $X$ and the
set of colored cones of $\mathbb{F}(X)$ such that for all two
$G$-orbits $\mathcal{O}_1$ and $\mathcal{O}_2$ in correspondence with
two colored cones $(\mathcal{C}_1,\mathcal{F}_1)$ and
$(\mathcal{C}_2,\mathcal{F}_2)$, we have
$\mathcal{O}_1\subset \overline{\mathcal{O}_2}$ if and only if
$(\mathcal{C}_2,\mathcal{F}_2)$ is a colored face of
$(\mathcal{C}_1,\mathcal{F}_1)$.

Let $\mathcal{O}$ be a $G$-orbit of $X$ associated to
$(\mathcal{C},\mathcal{F})\in\mathbb{F}(X)$. Then for every
$D\in\mathcal{F}$, we have $\overline{D}\supset\overline{\mathcal{O}}$.

One-codimensional $G$-orbits correspond to one-dimensional
colored cones of $\mathbb{F}(X)$ of the form $(\mathcal{C},\emptyset)$.And $G$-orbits of rank $n-1$ correspond to one-dimensional colored cones of $\mathbb{F}(X)$.
\item There exists a morphism between two $G/H$-embeddings $X$ and
  $X'$ if and only if for every colored cone
  $(\mathcal{C},\mathcal{F})$ of $\mathbb{F}(X)$ there exists a
  colored cone $(\mathcal{C'},\mathcal{F'})$ of $\mathbb{F}(X')$ such
  that $\mathcal{C}\subset\mathcal{C'}$ and
  $\mathcal{F}\subset\mathcal{F'}$. Moreover the morphism is proper if
  an only if the supports of $\mathbb{F}(X)$ and $\mathbb{F}(X')$ are
  the same.
\end{itemize}
\end{teo}

We will also need the following consequence of a characterisation of
Cartier divisors in spherical varieties, see \cite[Prop.3.1]{Br89}:

\begin{prop}\label{Qfactorial}
Let $X$ be a $\Qbb$-factorial spherical variety. Then each colored
cone of $\mathbb{F}(X)$ is simplicial. For any color $D$ of
$\mathbb{F}(X)$, $\rho(D)$ is in an edge of a cone of $\mathbb{F}(X)$
and if $\rho(D')$ is in the same edge, then $D=D'$ or $D'$ is not a
color of $\mathbb{F}(X)$.
\end{prop}

\subsection{Blow-ups of quasi-regular varieties}

We start with the following lemma:

\begin{lem}
\label{cone-ev} Let $X$ be a spherical $G$-variety of rank $n$. 
Then the following conditions are equivalent:
\begin{itemize}
\item[{\rm (QR)}] the irreducible components of $\partial X$ are of rank $n-1$ and
any $G$-orbit closure is the intersection of irreducible components
of $\partial X$ containing it;
\item[{\rm (VC)}] For any color $D$ of $\mathbb{F}(X)$, we have
  $\rho(D)\in\mathcal{V}$.
\end{itemize}
\end{lem}

\begin{proof}
Let us first remark, by Theorem \ref{Luna-Vust}, that the first part
of condition $(QR')$ is equivalent to say that the irreducible
components $Z_0,\dots,Z_r$ of $\partial X$ are the closures of the
$G$-orbits corresponding to the one dimensional colored cones of
$\mathbb{F}(X)$. This implies that, any $G$-orbit closure that is an
intersection of some $Z_i$'s, corresponds to a colored cone of
$\mathbb{F}(X)$ whose cone is generated by cones of one dimensional
colored cones of $\mathbb{F}(X)$.

Suppose that there exists a color $D$ of $\mathbb{F}(X)$ such that
$\rho(D)\not\in\mathcal{V}$. Let $(\mathcal{C},\mathcal{F})$ be a
colored cone of $\mathbb{F}(X)$ such that $D\in\mathcal{F}$. And let
$\mathcal{O}$ the $G$-orbit corresponding to
$(\mathcal{C},\mathcal{F})$. Then, by the preceding remark
$\overline{\mathcal{O}}$ is not the intersection of some $Z_i$'s.

Suppose now that condition (VC) holds. Let $\mathcal{O}$ be a
$G$-orbit and $(\mathcal{C},\mathcal{F})$ the corresponding colored
cone of $\mathbb{F}(X)$. By (VC), the cone $\mathcal{C}$ is
contained in $\mathcal{V}$, so that $(\mathcal{C})$ is generated by
the cones of its colored faces of dimension one. Then
$Z_0,\dots,Z_r$ corresponds to one dimensional colored cones of
$\mathbb{F}(X)$, so that they are of rank $n-1$. And any $G$-orbit
closure is the intersection of some $Z_i$'s.
\end{proof}

\begin{rem}
Remark that the equivalent conditions of the previous Lemma are
slightly weaker than quasi-regularity, the only difference being
that we do not assume the irreducible components of the boundary to
be smooth.
%
\end{rem}

\begin{lem}
\label{transverse} (\i) Let $X$ be a smooth projective spherical
variety such that $\partial X$ is the union of smooth irreducible
varieties $Z_0,\dots,Z_r$. Then the components $Z_0,\dots,Z_r$
intersect transversally.

(\i\i) In particular if $X$ is quasi-regular, then any closure of a
$G$-orbit is smooth.
\end{lem}

\begin{proof}
(\i) Let us use the local structure of spherical varieties
\cite{BLV}. Indeed, for any closed $G$-orbit $Y$ of $X$
there exists an affine open $B$-stable subvariety $X_0$ of $X$
intersecting $Y$ which is isomorphic to the product of the unipotent
radical $U_P$ of a parabolic subgroup $P$ and an affine $L$-stable
spherical subvariety $V$ of $X_0$ where $L$ is a Levi subgroup of
$P$. Moreover $V$ has a fixed point under $L$. Since $X$ is smooth,
$V$ is also smooth. By \cite{Lu73}, any smooth affine $L$-stable
spherical variety with a $L$-fixed point is an $L$-module. In
particular $V$ is an $L$-module.

Let $Z_i$ be an irreducible component of $\partial X$ containing
$Y$. Let $V_i$ be the intersection of $Z_i$ with $V$. The
intersection $Z_i\cap X_0$ is isomorphic to the product of $U_P$
with $V_i$. In particular $V_i$ is a proper $L$-stable smooth
irreducible subvariety of $V$. By the same argument, $V_i$ is a
sub-$L$-module of $V$. Writing the decomposition of $V$ into
irreducible submodules $V=\bigoplus_{k\in K}V(\chi_k)$, there exists
subsets $K_i$ of $K$ such that $$V_i=\bigoplus_{k\not\in
K_i}V(\chi_k).$$ Let us prove that the $K_i$, for
$i\in\{1,\dots,r\}$,  are disjoint. This will imply that
$V_0,\dots,V_r$, and then $Z_1,\dots,Z_r$, intersect transversally.
Suppose that there exists $i\neq j$ such that $K_i\cap
K_j\neq\emptyset$. Then $V_i$ and $V_j$ are in a same proper
sub-$L$-module of $V$. It implies that $Z_i$ and $Z_j$ are included
in a proper $P$-stable irreducible subvariety of $X$ (that is not
$G$-stable by maximality of $Z_i$ and $Z_j$). Since
$(G/H)\setminus(BH/H)$ is the union of the colors of $G/H$, we
proved that $Z_i$ and $Z_j$ are included in the closure of a color
of $G/H$. This is a contradiction by Theorem \ref{Luna-Vust} and
Proposition \ref{Qfactorial}.

(\i\i) Follows from (\i) and (QR2).
\end{proof}

\begin{prop}
\label{reg-quasi-reg} Let $X$ be a quasi-regular variety. Then there
exist quasi-regular varieties $X_0,\dots,X_r$ and morphisms
$\phi_i:X_i\to X_{i-1}$ for all $i\in\{1,\dots,r\}$ such that
$X_0=X$, $X_r$ is regular and $\phi_i$ is the blow-up of an
irreducible component of $\partial X_{i-1}$.
\end{prop}

\begin{proof}
Let $X$ be a quasi-regular variety. We proceed by induction on the
number of colors of $X$ (which is always finite). If $X$ has no color then
$X$ is regular and there is nothing to prove. Suppose that $X$ has a
color $D$. Since $X$ is smooth,
by Proposition \ref{Qfactorial}, $\rho(D)$ is in an edge of
$\mathbb{F}(X)$. This implies, by Theorem \ref{Luna-Vust}, that
$\overline{D}$ contains an irreducible component $Z$ of
$\partial X$. Let $\pi:\Xt\to X$ be the blow-up of $Z$ in $X$ and
let $E$ be the exceptional divisor. Because $X$ is quasi-regular, the
component $Z$ is smooth and $\Xt$ is a smooth spherical variety.
Remark that $D$ is not a color of $\Xt$ because $\overline{D}$ does
not contain $E$. Moreover if $D'$ is a color of $X$ different from
$D$, $\rho(D')$ is in another edge of $\mathbb{F}(X)$, so
$\overline{D'}$  contains a $G$-orbit of $X\backslash Z$ and then $D'$
is also a color of $\Xt$.

Now let us use the description of morphism between $G/H$-embeddings
given in Theorem \ref{Luna-Vust}, to conclude that $\mathbb{F}(\Xt)$
is obtained from $\mathbb{F}(X)$ by removing the color $D$ in all
colored cones that contain this color. Indeed, since $X$ is smooth,
the colored cones of $\mathbb{F}(X)$ are simplicial. Then, if
$\mathbb{F}(\Xt)$ is not obtained from $\mathbb{F}(X)$ by removing
some colors, there exists a one dimensional colored cone of $X$ not
in $\mathbb{F}(X)$. This gives a contradiction because $\partial X$
and $\partial \Xt$ have the same number of irreducible components so
that $\mathbb{F}(\Xt)$ and $\mathbb{F}(X)$ have the same number of
one dimensional colored cones.

The boundary of $\Xt$ is again the union of smooth irreducible
varieties. Indeed, these components are the strict transforms of
boundary components in $X$ and the exceptional divisor. Furthermore,
because of Lemma \ref{cone-ev} and the fact that our blow-up only
removes one color, the second condition for quasi-regularity is also
satisfied.
\end{proof}

\begin{rem}
Note that, if $X$ is a spherical variety, then there exists a
morphism $Y\to X$ which is a sequence of blow-ups of smooth
$G$-stable subvarieties such that $Y$ is regular. This result
follows from a general result of \cite{RY} on rational
$G$-equivariant morphism (see \cite[Proof of Corollary
4.4.2]{brion}). But in general, the smooth $G$-stable subvarieties
that are blown-up, are not irreducible components of $\partial X$.
For example, if Condition (VC) is not satisfied, to remove the color
whose image in $N$ is outside the valuation cone $\mathcal{V}$, we
must blow-up $X$ along a $G$-stable subvariety that is not an
irreducible component of $\partial X$.
\end{rem}

To prove Theorem \ref{eclatement}, we need to study the behaviour of
the action sheaf under successive blow-ups with smooth centers. We
do this in the next subsection.

\subsection{Action sheaf and blow-ups}

Let $X$ be a smooth variety and $Y$ be any smooth subvariety in $X$.
Let us denote by $N_Y$ the normal sheaf $Y$ in $X$. It this
situation, it is a vector bundle on $Y$ and there is a natural
surjective morphism $T_X\to N_Y$ (here by abuse of notation we still
denote by $N_Y$ the push-forward of the normal bundle by the
inclusion of $Y$ in $X$). The action sheaf $T_{X,Y}$ of $Y$ in $X$
is the kernel of this map. In symbols, we have an exact sequence:
$$0\to T_{X,Y}\to T_X\to N_Y\to 0.$$
If furthermore $Y$ is of codimension 1 (i.e. a
Cartier divisor), then it is not difficult to see that the action sheaf
is locally free (see for example \cite{BB96}).

\vskip 0.2 cm

Consider the blow-up $\pi:\Xt\to X$ of $Y$ in $X$. Let us denote by
$E$ the exceptional divisor. In this subsection, we want to relate
the action sheaf $T_{X,Y}$ with the push forward $\pi_* T_{\Xt,E}$
of the action sheaf of $E$ in $\Xt$. We shall prove the following

\begin{lem}
 With the notations above, we have the equalities
\begin{itemize}
\item $\pi_*(T_{\Xt,E})=T_{X,Y}$ and
\item $R^i\pi_*(T_{\Xt,E})=0$ for all $i>0$.
\end{itemize}
\end{lem}

\begin{proof}
Recall the definition of the tautological quotient bundle $Q$ on the
exceptional divisor $E$ given by $Q=\pi^*N_Y/\cO_E(-1)$. The
following exact sequence holds (see for example \cite[Page
299]{Fu98}):
$$0\to T_{\Xt}\to\pi^* T_X\to Q\to 0.$$
In particular, the composition of the differential of $\pi$ given by
$T_\Xt\to\pi^* T_X$ and the map $\pi^* T_X\to\pi^*N_Y$ to the normal
bundle of $Y$ factors through $\cO_E(-1)$ the normal bundle of $E$.
We thus have a commutative diagram
$$\xymatrix{&&0\ar[d]&0\ar[d]&\\
0\ar[r]& T_{\Xt,E}\ar[r]& T_{\Xt}\ar[r]\ar[d]&
\cO_E(-1)\ar[r]\ar[d]& 0\\
&\pi^*T_{X,Y}\ar[r]&\pi^*T_X\ar[r]\ar[d]&\pi^*N_Y\ar[r]\ar[d]& 0\\
&&Q\ar[d]\ar@{=}[r]&Q\ar[d]&\\
&&0&0&}.$$
proving that $T_{\Xt,E}$ is the kernel of the map
$\pi^*T_X\to\pi^*N_Y$. Pushing forward the associated exact sequence
using that $\pi_*\cO_\Xt=\cO_X$ and the fact that the pushed forward
map $T_X\to N_Y$ is surjective, the result follows.
\end{proof}

\begin{rem}
\label{calcul_local}
 Remark that the sheaves $T_{\Xt,E}$ and $\pi^*T_{X,Y}$ do not
coincide on $\Xt$. Indeed, the sheaf $T_{\Xt,E}$ is locally free
(because $E$ is a divisor) whereas the sheaf $\pi^*T_{X,Y}$ has a
non trivial torsion part. To see this we may compute in local
coordinates in an {\'e}tale neigbourhood of a point in $Y$ (or
equivalently an open neigbourhood for the usual topology). In such
an open subset $U$, the ring of $X$ can be chosen to be
$k[U]=k[x_1,\cdots,x_n]$ where the ideal of $Y$ is
$(x_1,\cdots,x_p)$. We may choose an open subset $V$ above $U$ such
that the coordinated ring of $\Xt$ is given by
$k[V]=k[y_1,\cdots,y_{p-1},x_p,\cdots,x_n]$ and the map $V\to U$
gives rise to a morphism $k[U]\to k[V]$ given by $x_i\mapsto x_i$
for $i\geq p$ and $x_i\mapsto y_ix_p$ for $i\leq p-1$.

In these coordinates, the sheaf $T_X$ is generated by the vectors
$({\frac{\partial}{\partial x_i}})_{i\in[1,n]}$ while the sheaf
$T_{X,Y}$ is generated by the vectors
$m_{i,j}=x_i\frac{\partial}{\partial x_j}$ for $i$ and $j$ in $[1,p]$
and the vectors $\frac{\partial}{\partial x_i}$ for $i>p$. These
vectors satisfy the relations
$$x_km_{i,j}=x_im_{k,j}$$
for all $i$, $j$ and $k$ in $[1,p]$. In particular if we pull this
sheaf back to $V$ we may consider the elements $m_{i,j}-y_im_{p,j}$
multiplying by $x_p$ we get
$$x_p(m_{i,j}-y_im_{p,j})=x_pm_{i,j}-x_py_im_{p,j}=
x_pm_{i,j}-x_im_{p,j}=0$$
and the element $m_{i,j}-y_im_{p,j}$ is a torsion element.

The quotient of the sheaf $\pi^*T_{X,Y}$ by its torsion part is the
image in $\pi^*T_X$ of the sheaf $\pi^*T_{X,Y}$. It is also easy in
these coordinates to prove that the images of $\pi^*T_{X,Y}$ and
$T_{\Xt,E}$ in $\pi^*T_X$ coincide reproving the previous
lemma. Indeed, let us first compute the image of $\pi^*T_{X,E}$ in
$\pi^*T_X$. It is generated by the vectors
$x_i\frac{\partial}{\partial x_j}$ for $i$ and $j$ in $[1,p]$ and the
vectors $\frac{\partial}{\partial x_k}$ for $k>p$. Because
$x_i=y_ix_p$ for $i<p$ we get the following set of generators for this
image:
$$x_p\frac{\partial}{\partial x_j},\ \textrm{for $j\leq p$}\ \textrm{  and
}\ \frac{\partial}{\partial x_k} \textrm{ for $k>p$}.$$
In particular we see that it is a locally free subsheaf of
$\pi^*T_X$. To determine the image of $T_{\Xt,E}$, let us first recall
the map $d\pi:T_\Xt\to\pi^*T_X$. The sheaf $T_\Xt$ is locally free
generated by the vectors $(\frac{\partial}{\partial
  y_i})_{i\in[1,p-1]}$ and $(\frac{\partial}{\partial
  x_k})_{k\in[p,n]}$. The differential of $\pi$ is defined as follows:
$$\frac{\partial}{\partial y_i}\mapsto x_p\frac{\partial}{\partial
  x_i} \textrm{for $i<p$, } \frac{\partial}{\partial x_p}\mapsto
\frac{\partial}{\partial x_p}+
\sum_{i=1}^{p-1}y_i\frac{\partial}{\partial x_i} \textrm{ and }
\frac{\partial}{\partial x_k}\mapsto \frac{\partial}{\partial x_k}
\textrm{for $k>p$}.$$
The action sheaf $T_{\Xt,E}$ is generated by the vectors
$(\frac{\partial}{\partial y_i})_{i\in[1,p-1]}$,
$x_p\frac{\partial}{\partial x_p}$ and $(\frac{\partial}{\partial
  x_k})_{k\in[p+1,n]}$. Its image is thus the same as the image of
$\pi^*T_{X,Y}$.
\end{rem}

More generally, we want to study the action sheaf of a subvariety $Z$ of $X$
containing $Y$ as an irreducible component and such that the other
components of $Z$ meet $Y$ transversally. More precisely, let us
write the decomposition into irreducible components of $Z$ as follows:
$$Z=\bigcup_{i=0}^rZ_i$$
with $Z_0=Y$ and assume that all these components intersect
transversally. We shall prove the following:

\begin{prop}
\label{eclt}
  Let $Z$ be a reduced subvariety of $X$ whose irreducible components
are smooth
  and containing $Y$ as an irreducible component. Assume furthermore
  that the irreducible components of $Z$ meet transversally. Let us
  denote by $Z'$ union of the irreducible components in $Z$ different
  from $Y$ and by $\Zt'$ the strict transform of $Z'$ by $\pi:\Xt\to
  X$ the blowing-up of $X$ along $Y$. Finally, let us denote by $\Zt$
  the union of $\Zt'$ and $E$ the exceptional divisor. Then we have
  the equalities
\begin{itemize}
\item $\pi_* (T_{\Xt,\Zt})=T_{X,Z}$ and
\item $R^i\pi_*(T_{\Xt,\Zt})=0$ for all $i>0$.
\end{itemize}
\end{prop}

\begin{proof}
  To prove this result, we only need --- as in the case of the
  previous lemma where $Z=Y$ --- to prove that the images of
  $\pi^*T_{X,Z}$ and of $T_{\Xt,\Zt}$ in $\pi^*T_X$ coincide. Indeed,
  we get in that case an exact sequence
$$0\to T_{\Xt,\Zt}\to\pi^*T_X\to \pi^*N_Z\to 0$$
and the result follows by push forward. To prove the equality of
these two images, we can use the local coordinate description one
more time. Indeed, because $Y$ and the other components of $Z$ meet
transversally, the equations of the other components
$Z_1,\cdots,Z_r$ of $Z$ in the coordinate ring
$k[U]=k[x_1,\cdots,x_n]$ may be chosen to be
$x_{p+1}=\cdots=x_{p_1}=0$ for $Z_1$ and more generally
$x_{p_{i-1}+1}=\cdots=x_{p_i}=0$ for $Z_i$ where
$p=p_0<p_1<\cdots<p_r\leq n$. In particular the computation made in
Remark \ref{calcul_local} do not change and the result follows.
\end{proof}

\begin{preu}
To conclude the proof of Theorem \ref{eclatement}, we use
Proposition \ref{reg-quasi-reg} and Lemma \ref{transverse} to prove
that we are in the situation of Proposition \ref{eclt}.

We conclude by using Theorem \ref{BBK} and the classical fact that
for a morphism $f:X\to Y$ and a sheaf $\mathcal{F}$ on $X$, if we have
$R^if_*\mathcal{F}=0$ for $i>0$, then
$H^i(X,\mathcal{F})=H^i(Y,f_*\mathcal{F})$.
\end{preu}

Now let us give examples of quasi-regular varieties. The following
result implies that the result of Theorem \ref{eclatement} applies
to smooth horospherical varieties and to smooth spherical varieties
of rank one:

\begin{lem}\label{horospherical}
(\i) Smooth projective horospherical varieties are quasi-regular.

(\i\i) Smooth projective spherical varieties of rank~1 are
quasi-regular.
\end{lem}

\begin{proof}
(\i) Note that any $G$-stable subvariety of a smooth horospherical
variety is smooth \cite[Chap.2]{Pa06} then (QR1) is
satisfied. Moreover, if the homogeneous space $G/H$ is horospherical,
the valuation cone is the vector space $N_\Qbb$ \cite[Cor.7.2]{Kn91},
so that (VC) is automatically satisfied and we conclude by Lemma
\ref{cone-ev}.

(\i\i) Let $X$ a smooth projective spherical variety of rank~1. We
can assume that $X$ is not horospherical. Then $X$ has two orbits,
so that $\partial X$ is a flag variety. Conditions (QR1) and (QR2)
are clearly satisfied.
\end{proof}

\begin{rem} These results cannot be extended in the same way to
spherical varieties of rank more than~1. Indeed, a $G$-stable
subvariety of a smooth spherical variety is not necessarily smooth
(see \cite{Br94}).
\end{rem}

\section{Applications to spherical varieties of rank~1}
\subsection{Proof of Theorem \ref{deformations}}

If $G$ is a connected reductive algebraic group, we denote by $B$ a
Borel subgroup of $G$ containing a maximal torus $T$. We denote by
$\varpi_i$ the fundamental weights of $(G,B,T)$ with the notation of
\cite{Bo75}, and by $P(\varpi_i)$ the corresponding maximal parabolic
subgroup containing $B$. Let $\rho$ be the sum of the fundamental
weights.

Let $P$ be a parabolic subgroup of $G$ and $V$ a $P$-module.  Then
the homogeneous vector bundle $G\times^PV$ over $G/P$ is the
quotient of the product $G\times V$ by the equivalence
relation~$\sim$ defined by $$\forall g\in G,\,\forall p\in P,
\forall v\in V,\quad (g,v)\sim(gp^{-1},p.v).$$ To compute the
cohomology of such vector bundles on flag varieties, we will use the
Borel-Weil Theorem (see \cite[Chapter 4.3]{ak}):

\begin{teo}[Borel-Weil]
Let $V$ be an irreducible $P$-module of highest weight $\chi$.
Denote by $\mathcal{V}(\chi)$ the vector bundle $G\times^PV$ over
$G/P$ and by $w_0^P(\chi)$ the lowest weight of $V$. We have the
following alternative:
\begin{itemize}
\item  If there exists a root $\alpha$ 
with $\langle w_0^P(\chi)-\rho,\alpha^\vee\rangle=0$, then, for
any $i\geq 0$, $H^i(G/P,\mathcal{V}(\chi))=0$.
\item Otherwise, there exists an element $w$ of the Weyl group 
with $\langle w(w_0^P(\chi)-\rho),\alpha^\vee\rangle<0$ for all
positive roots $\alpha$. Denote by $l(w)$ the length of $w$. Then we
have $H^i(G/P,\mathcal{V}(\chi))=0$ for $i\neq l(w)$ and
$H^{l(w)}(G/P, \mathcal{V}(\chi))$ is the $G$-module of highest
weight $-w(w_0^P(\chi)-\rho)-\rho$.
\end{itemize}
\end{teo}

In this section we shall freely use the results of \cite{Pa08}. Let
$X$ be one of the horospherical varieties satisfying the condition
$(\dag)$ in the introduction (recall that there varieties are
$X^1(m)$, $X^2$, $X^3(m,i)$, $X^4$ and $X^5$).  Denote by $P_Y$ and
$P_Z$ the maximal parabolic subgroups of $G$ containing $B$ such
that the closed orbits $Y$ and $Z$ of $X$ are respectively
isomorphic to $G/P_Y$ and $G/P_Z$. According to \cite[Section
1.4]{Pa08}, there exists a character $\chi$ of $P_Y\cap P_Z$ such
that the total spaces of the normal bundles $N_{Y}$ and $N_{Z}$ are
respectively $G\times^{P_Y}V(\chi)$ and $G\times^{P_Z}V(-\chi)$.
Then, to compute the cohomology of $N_{Y}$ and $N_{Z}$ applying
Borel-Weil Theorem, we only need to compute the lowest weights
$w_0^Y(\chi)$  of the $P_Y$-module $V(\chi)$ and $w_0^Z(-\chi)$ of
the $P_Z$-module $V(-\chi)$. Remark that, if $P$ is a parabolic
subgroup of $G$, $w_0^P$ is the longest element of the Weyl group
fixing the characters of $P$.
We summarize in the following table, the value of these data in each case.
\vskip 0.1 cm
\begin{center}
\begin{tabular}{|c|c|c c|c|c c|}
\hline
 $X$ & Type of $G$ & $P_Y$ & $P_Z$ & $\chi$ &  $w_0^Y(\chi)$ & $w_0^Z(-\chi)$\\
\hline
\hline
$X^1(m)$ & $B_m$ & $P(\varpi_{m-1})$ & $P(\varpi_m)$ &
$\varpi_m-\varpi_{m-1}$ & $-\varpi_m$ & $-\varpi_1+\varpi_m$\\
\hline
$X^2$ & $B_3$ & $P(\varpi_1)$ & $P(\varpi_3)$ & $-\varpi_3$ &
$\varpi_3-\varpi_1$ & $-\varpi_2+\varpi_3$\\
\hline
$X^3(m,i)$ & $C_m$ & $P(\varpi_i)$ & $P(\varpi_{i+1})$ &
$\varpi_{i+1}-\varpi_i$ & $-\varpi_{i+1}+\varpi_i$ & $-\varpi_1$\\
\hline
$X^4$ & $F_4$ & $P(\varpi_2)$ & $P(\varpi_3)$  & $\varpi_3-\varpi_2$ &
$-\varpi_4$ & $-\varpi_1+\varpi_3$\\
\hline
$X^5$ & $G_2$ & $P(\varpi_1)$ & $P(\varpi_2)$ & $\varpi_2-\varpi_1$ &
$-\varpi_2+2\varpi_1$ & $-\varpi_1$\\
\hline
\end{tabular}
\end{center}
\vskip 0.1 cm
In each case, by the Borel-Weil Theorem, one of the normal bundle
has cohomology in degree~0 (and only in degree~0). This was already
used in \cite{Pa08}. And the other normal bundle has no cohomology
except in the last case where $N_{Y}$ has cohomology in degree~1
(and only in degree~1). Indeed, when $X=X^5$, we have
$H^1(X,N_{Y})=\Cbb.$

When $X$ is one of the two spherical varieties $\mathbb{X}_1$ and
$\mathbb{X}_2$, denote by $P_Y$ the parabolic subgroup of $G$
containing $B$ such that the unique closed orbit $Y$ is isomorphic to
$G/P_Y$.  Let $\chi$ be the character of $P_Y$ such that the total
space of the normal bundle $N_{Y}$ is  $G\times^{P_Y}V(\chi)$ and let
$w_0^Y(\chi)$ be the lowest weight of the $P_Y$-module $V(\chi)$. Then
we have the following table.
\vskip 0.1 cm
\begin{center}
\begin{tabular}{|c|c|c|c|c|}
\hline
 $X$ & $G$ & $P_Y$ & $\chi$ &  $w_0^Y(\chi)$\\
\hline
\hline
$\mathbb{X}_1$ & $F_4$ & $P(\varpi_3)$ & $\varpi_1-\varpi_3$ &
$-\varpi_2+\varpi_3$\\
\hline
$\mathbb{X}_2$ & $G_2\times\operatorname{PSL(2)}$ & $P(\varpi_1)\cap
P(\varpi_0)$ & $\varpi_2-2\varpi_1-2\varpi_0$ &
$-\varpi_2+\varpi_1-2\varpi_0$\\
\hline
\end{tabular}
\end{center}
\vskip 0.1 cm In both cases, the Borel-Weil Theorem implies the
vanishing of the cohomology of the normal bundle of $Y$ in $X$.

Then Theorem \ref{deformations} is a corollary of Theorem
\ref{eclatement} using the long exact sequence defining the action
sheaf:
$$0\lra S_X\lra T_X\lra N_{\partial X/X}\lra 0$$ where $\partial X=Y\cup
Z$ if $X$ is horospherical and $Y$ in the other cases.

\begin{rem}
This result gives an example of an horospherical variety such that
$H^1(X,T_{X,Y})\neq 0$ with $Y$ is a $G$-stable subvariety of $X$
(different from $\partial X$). Indeed, let $X$ be the smooth
projective non-homogeneous horospherical $G_2$-variety of Picard
number one (\emph{i.e.} the variety $X^5$). We have the following
short exact sequence
$$0\lra S_X\lra T_{X,Y}\lra N_{Z}\lra 0.$$
This implies the equalities $H^1(X,T_{X,Y})=H^1(X,N_{Z})=\Cbb$.
\end{rem}

\subsection{Explicit deformation of the horospherical $G_2$-variety $X^5$}

In the unique case where $X$ is not locally rigid, we prove that the
local deformation comes from a global deformation of $X$. We describe
explicitly the unique deformation of $X$ in the following proposition.

\begin{prop}
\label{deformation} The variety $X^5$ has a deformation to the
orthogonal grassmannian $\operatorname{Gr}_q(2,7)$.
\end{prop}

\begin{proof} Let us denote by
  $(z_0,z_1,z_2,z_3,z_{-1},z_{-2},z_{-3})$ the basis of the imaginary
  octonions $\Im(\Obb)$ as in \cite[Section.2.3]{Pa08}. Then we have
  $$X^5=\overline{G_2.[z_3+z_1\wedge z_3]}\subset \Pbb(V(\varpi_1)\oplus
  V(\varpi_2)).$$
Now for all $t\in\Cbb^*$, define $x_t:=z_3+t(z_0\wedge z_3+z_1\wedge
z_2)+z_1\wedge z_3\in V(\varpi_1)\oplus V(\varpi_2)$
and set
$$X^5_t:=\overline{G_2.[x_t]}\subset\Pbb(V(\varpi_1)\oplus V(\varpi_2)).$$
The limit of the varieties $X_t$ when $t$ goes to 0 is $X^5$. We
conclude the proof thanks to the following lemma.
\end{proof}

\begin{lem}
For all $t\in\Cbb^*$, $X_t$ is isomorphic to the orthogonal
grassmannian $\operatorname{Gr}_q(2,7)$.
\end{lem}

\begin{proof}
For $\tau$ in $\mathbb{C}$, let $\phi_\tau\in G_2$ be the
automorphism of the $\Obb$ sending $(1,z_0,z_1,z_2,z_3, z_{-1},
z_{-2},z_{-3})$ to $(1,z_0+\tau
z_1,z_1,z_2+\frac{\tau}{2}z_3,z_3,z_{-1}+\tau
z_0+\frac{\tau^2}{2}z_1,z_{-2},z_{-3}-\tau z_{-2})$. Then we have
that the element $\phi_{-2/3t}(x_t)=z_3+t(z_0\wedge z_3+z_1\wedge
z_2)$ is in the open orbit of $X_t$. This proves the lemma, by
\cite[Prop 2.34]{Pa08}.
\end{proof}

\subsection{Automorphisms}

%
%
%
%
%
%

In \cite{Pa08}, a case by case analysis was necessary to prove that
certain spherical varieties of rank~1 are homogeneous. The following
result gives a more uniform way to check the homogeneity these
varieties.

\begin{prop}
\label{aut} Let $X$ be a smooth projective spherical variety of
rank~1. Then $\Aut(X)$ acts transitively on $X$ if and only if for
all closed $G$-orbit $Y$ of $X$, the normal bundle $N_{Y/X}$ has a
non-zero section.
\end{prop}

\begin{proof}
Let $X$ be a smooth spherical variety of rank~1. Recall that, in that
case, $\partial X$ is one closed orbit or the union of two-closed
orbits.

Suppose that there exists a closed $G$-orbit $Y$ of $X$ such that
$H^0(X,N_{Y})=0$ then we have the equality $H^0(X,T_{X,Y})=H^0(X,T_X)$
and $\Aut(X)$ stabilises $Y$.

Suppose now that, for all closed $G$-orbit $Y$ of $X$, the normal
bundle $N_{Y}$ has a non-zero section. This implies, by the Borel-Weil
Theorem, that for all closed $G$-orbit $Y$ of $X$, we have
$H^1(X,N_{Y})=0$. We thus have the following exact sequence $$0\lra
H^0(X,T_{X,Y}) \lra H^0(X,T_X)\lra H^0(Y,N_{Y})\lra 0.$$ Indeed, if
$X$ has only one closed $G$-orbit, then $H^1(X,T_{X,Y})=0$ because
$T_{X,Y}=S_X$. If $X$ has two closed $G$-orbits $Y$ and $Z$, then we
have the following exact sequence
$$0\to S_X \to T_{X,Y}\to N_{Z}\to 0$$
and we obtain $H^1(X,T_{X,Y})=0$ from Theorem \ref{eclatement}.
This implies that $Y$ cannot be stabilised by $\Aut(X)$ and the result
follows.
\end{proof}

\bigskip\noindent

\medskip\noindent
Boris {\sc Pasquier}, Hausdorff Center for Mathematics,
Universit{\"a}t Bonn, Landwirtschaftskammer (Neubau)
Endenicher Allee 60, 53115 Bonn, Germany.

\noindent {\it email}: \texttt{boris.pasquier@hcm.uni-bonn.de}

\medskip\noindent
Nicolas {\sc Perrin}, Hausdorff Center for Mathematics,
Universit{\"a}t Bonn, Landwirtschaftskammer (Neubau) Endenicher
Allee 60, 53115 Bonn, Germany and Institut de Math{\'e}matiques de
Jussieu, Universit{\'e} Pierre et Marie Curie, Case 247, 4 place
Jussieu, 75252 PARIS Cedex 05, France.

\noindent {\it email}: \texttt{nicolas.perrin@hcm.uni-bonn.de}

\end{document}